# Directional Congestion in Data Envelopment Analysis

Guo-liang Yang

Institute of Policy and Management, Chinese Academy of Sciences, Beijing 100190, China

**Abstract.** First, this paper proposes the definition of directional congestion in certain input and output directions in the framework of data envelopment analysis (DEA). Second, two methods from different viewpoints are also proposed to estimate the directional congestion in a DEA framework. Third, we address the relations among directional congestion and classic congestion and strong (weak) congestion. Finally, we present a case study investigating the analysis of the research institutes in the Chinese Academy of Sciences (CAS) to demonstrate the applicability and usefulness of the methods developed in this paper.

## 1 Introduction

Congestion is often involved in the real practice of analyzing the returns to scale (RTS) of decision making units (DMUs), which describes the case whereby the decrease of one (or some) inputs will cause the maximum possible increase of one (or some) outputs without worsening any other input or output (Cooper et al., 2004). Essentially, the congestion effect reflects the problem of excessive inputs (Wei and Yan, 2004). We take as an example the production process with a single input and single output, which is exhibited in Figure 1.

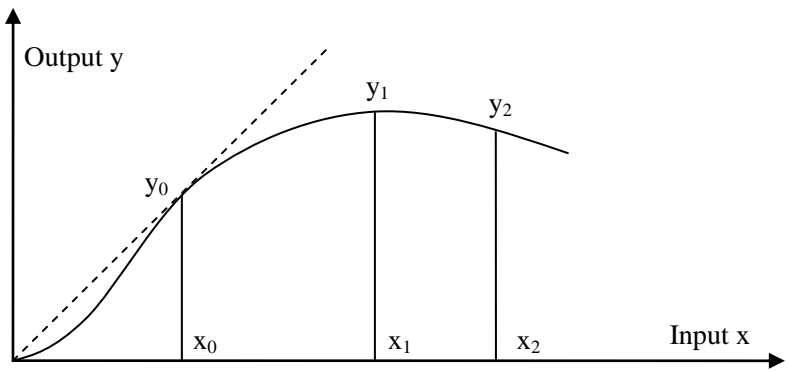

**Figure 1: Production function with a single input and single output.**

In Figure 1, the horizontal axis represents the input x, and the vertical axis represents the output y. The curve is a single-peak function, which represents the maximum possible output generated from a certain number of inputs.

Because the point on the curve represents the maximum possible output by a certain level of input, we can see the point on the curve is an efficient production process. Congestion represents the falling portion of the curve, i.e., the (y1-y2) segment. Thus, the definition of congestion in economics is as follows (Cooper et al., 2001a; Brockett et al., 2004):

**Definition 1 (Congestion):** Evidence of congestion is present when the reductions in one or more inputs can be associated with the increases in one or more outputs without worsening any other input or output. Proceeding in reverse, congestion occurs when the increases in one or more inputs can be associated with the decreases in one or more outputs without improving any other input or output.

Cooper et al. (2004) noted that the research on congestion in classic economics is insufficient, partly due to the comments put forward by the Nobel laureate Stigler (1976) on the "X-Efficiency" proposed by Leibenstein (1966, 1976) questioning whether congestion should be a research topic in economics. However, Färe and Svensson (1980) redefined the concept of congestion. Subsequently, Färe and Grosskopf (1983) proposed the operational concept of congestion. Färe et al. (1985) first proposed the corresponding DEA (data envelopment analysis) model to explore the congestion effect. Subsequently, Cooper et al. (1996) proposed an alternative DEA approach to investigate the congestion effect. Cooper et al. (2001b) compared the similarities and differences between these two models. Based on the weak disposal assumption, Wei and Yan (2004) and Tone and Sahoo (2004) rebuilt the production possibility set (PPS) and the corresponding DEA model to determine the congestion effect of the DMUs. Kao (2010) investigated three types of models prevalent in the DEA literatures for measuring the congestion effect and utilized the model proposed by Wei and Yan (2004) on Taiwan forests. Based on the measurement and reorganization, he gave some interesting findings on how to alleviate the congestion in Taiwan forests. Brockett et al. (1998) used the model proposed by Cooper et al. (1996) to identify congestion inputs in Chinese industries. Shortly afterward, Cooper et al. (2001a) examined the problem of inefficiency in the Chinese automobile and textile industries. Using the annual data for the period 1981 to 1997, they found the evidence of congestion in both industries. Jahanshahloo and Khodabakhshi (2004) introduced a suitable combination of inputs for improving outputs with determining input congestion. Sueyoshi and Sekitani (2009) explored how to deal with the occurrence of multiple solutions in the DEA-based congestion measurement. Noura et al. (2010) presented a new method for measuring congestion. Khoveyni et al. (2013) proposed a slack-based DEA approach to recognize congestion (strong and weak) for the target DMUs with non-negative inputs and outputs. To date, there have been three main approaches on the measurement of congestion in the DEA framework. Please see

Section 2 for details.

The existing studies on the RTS or congestion measurement in the DEA models are all based on the definition of the RTS in the DEA framework based on the PPS made by Banker (1984), which extended the application area of the DEA from a relative efficiency evaluation to a RTS measurement. However, Podinovski and Førsund (2010) and Atici and Podinovski (2012) indicated that the derivate in the definition in Banker (1984) may not always exist, so they replaced the classical derivative by the directional derivatives. However, for the point on the production function that is not differentiable, there always exist directional derivatives at this point, i.e., we can define the left- and right-hand directional SE, respectively. According to their value, we can determine the increasing (or the constant or decreasing) RTS on the left-hand (or right-hand) of the DMUs in a certain direction of inputs and outputs. Subsequently, Yang (2012) and Yang et al. (2014) extended the classic definition of RTS in economics to directional RTS and introduced the definition into the DEA framework. He extended the classic definitions of SE or RTS in the DEA proposed by Banker (1984) and redefined the left- and right-hand direction SE and the directional RTS based on the PPS.

There exist three concepts describing congestion: congestion, strong congestion and weak congestion, which are defined as **Definition 1, 2 and 3**, respectively (see Section 2 for **Definition 2 and 3**). In existing congestion measurement, strong congestion and weak congestion are defined to describe the case of either an activity that uses less inputs (all inputs decrease proportionally) to produce more outputs (all outputs increase proportionally) or an activity that uses less resources in one or more inputs to make more products in one or more outputs. We can observe that strong congestion is a special case of weak congestion. Sueyoshi and Sekitani (2009) noted that there are the following relationships among the three concepts (congestion, weak congestion and strong congestion): *"(a) If a DMU belongs to strong congestion, then the DMU belongs to congestion. (b) If a DMU is strongly efficient with respect to $P_{convex}(X,Y)$ and it belongs to congestion, then the DMU belongs to weak congestion"*.

There are two main problems in the existing definitions of strong and weak congestions. First, from their definitions we observe that strong congestion and weak congestion describe the case of input decrease, which is the former part of **Definition 1**. These two definitions fail to tell us what will happen when increases in one or more inputs are associated with decreases in one or more outputs – without improving any other input or output, which is the latter part of **Definition 1**. Second, if weak congestion occurs, we do not know the precise direction in which the congestion occurs. We can only know the slacks of corresponding models instead of the directions of the inputs or outputs (i.e., the proportions of them). This means we do not know whether congestion will occur if decision-makers (DMs) decide to increase (or decrease) the inputs disproportionally. To this end, this paper aims to introduce the definition of directional congestion to describe whether congestion occurs in certain input and output directions and, subsequently, to propose the corresponding quantitative methods to measure the extent of directional congestion in the DEA framework.

The rest of this paper is organized as follows. Section 2 summarizes the three main approaches on congestion measurement. Section 3 proposes the definitions of directional congestion. Two approaches, the Finite Difference Method (FDM) and the Upper and Lower Bounds Method (ULBM), are proposed in Section 4 to estimate the directional congestion. Furthermore, in this section we discuss certain characteristics of directional congestion compared with classic congestion. Section 5 provides an illustrative example of the in-depth analysis of 16 basic CAS research institutes in 2009. The last section presents the conclusions.

## 2 Primary Approaches on Congestion Measurement

Suppose that there are n DMUs (j=1,…,n) to be analyzed. Each $DMU_j$ has m inputs and s outputs, which are denoted by $x_{ij}$ (i=1,…, m) and $y_{rj}$ (r=1,…,s), respectively. All the values of these indicators are not negative, and at least one is larger than zero. $DMU_0$ is the DMU to be evaluated and denoted as DMU $\left( X_0, Y_0 \right)$.

Currently, there are three main ideas on congestion research in the DEA framework. The first step of these three approaches is to measure the efficiencies of the DMUs using the BCC model proposed by Banker et al. (1984). Cooper et al. (2000) noted that the input-oriented BCC model may produce erroneous results when measuring the extent of congestion, so we use the output-based BCC model as the first step of the congestion analysis, which can be used to obtain the relative efficiency of $DMU_0$, i.e., the BCC-efficiency, and is listed as follows.

$$\hat{\theta} = \max_{\lambda_j, \theta, s_i^-, s_r^+} \left\{ \theta - \varepsilon \left( \sum_{i=1}^{m} s_i^- + \sum_{r=1}^{s} s_r^+ \right) \right\}$$

$$s.t. \begin{cases} \sum_{j=1}^{n} \lambda_j x_{ij} + s_i^- = x_{i0}, i = 1, \ldots, m \\ \sum_{j=1}^{n} \lambda_j y_{rj} - s_r^+ = \theta y_{r0}, r = 1, \ldots, s \\ \sum_{j=1}^{n} \lambda_j = 1 \\ s_i^-, s_r^+, \lambda_j \geq 0, i = 1, \ldots, m; r = 1, \ldots, s; j = 1, \ldots, n, \theta \text{ free} \end{cases} \quad (1)$$

where $\varepsilon$ is a non-Archimedean construct to ensure strongly efficient solutions, and variables $s_i^-, s_r^+$ are slacks.

**The first approach** can be called the FGL model (Färe et al., 1985). We know that strong disposability is implied in Model (1), i.e., if the inputs are increased, the outputs will not be reduced. Färe et al. (1985) proposed the assumption of weak disposability and constructed the model to estimate congestion using the following model:

$$\max_{\lambda_j, \beta, \tau} \beta$$

$$s.t. \begin{cases} \sum_{j=1}^{n} \lambda_j x_{ij} = \tau x_{i0}, i = 1, \ldots, m \\ \sum_{j=1}^{n} \lambda_j y_{rj} \geq \beta y_{r0}, r = 1, \ldots, s \\ \sum_{j=1}^{n} \lambda_j = 1 \\ 1 \geq \tau \geq 0, \lambda_j \geq 0, j = 1, \ldots, n, \beta \text{ free} \end{cases} \quad (2)$$

Let $\theta^*$ and $\beta^*$ represent the optimal objective value of Model (1) and Model (2), respectively. Byrnes et al. (1988) defined $1/\beta^*$ as pure technical efficiency (PTE)[1] and $\beta^*/\theta^*$ as the measurement to gauge the extent of congestion, which denotes the ratio between the relative efficiency under the assumption of strong disposal and that under the assumption of weak disposal. Thus, we have the following equation:

BCC-efficiency ($1/\theta^*$) = PTE ($1/\beta^*$) × Congestion ($\beta^*/\theta^*$)

It should be noted that Model (1) is equivalent to Model (2) when there is only one input, i.e., $\beta^*/\theta^* = 1$. In this case, there is no congestion in the FGL approach.

**The second approach** is the CTT model (Cooper et al., 1996). The CTT model is slack-based and its approach is to estimate the congestion effect through two stages. First, we determine the reference target of DMU$_0$ on the efficient frontier produced from Model (1) using the following formula:

$$\begin{cases} \hat{x}_{i0} = x_{i0} - s_i^{-*}, i=1,...,m, \\ \hat{y}_{r0} = \theta^* y_{r0} + s_r^{+*}, r=1,...,s, \end{cases}$$

where $(\theta^*, s_i^{-*}, s_r^{-*})$ is the optimal solution of Model (1). DMU$_0$ is BCC-efficient if and only if $\hat{x}_{i0} = x_{i0}$ and $\hat{y}_{r0} = y_{r0}, \forall i, r$, where $\forall$ means "for all".

---

[1] Note: The pure technical efficiency defined in Byrnes et al. (1984) is different from the efficiency measured in the BCC model, which is denoted by the BCC-efficiency.

To determine the extent of the input congestion, the second stage in the CTT model is to use the following Model (3) to calculate the maximal inputs that could be decreased.

$$\max_{\lambda_j, \delta_i^-} \sum_{i=1}^m \delta_i^-$$
$$s.t. \begin{cases} \sum_{j=1}^n \lambda_j x_{ij} - \delta_i^- = \hat{x}_{i0}, i=1,...,m \\ \sum_{j=1}^n \lambda_j y_{rj} = \hat{y}_{r0}, r=1,...,s \\ \sum_{j=1}^n \lambda_j = 1 \\ \delta_i^- \leq s_i^{-*} \\ \delta_i^-, \lambda_j \geq 0, i=1,...,m; j=1,...,n \end{cases} \quad (3)$$

where $\delta_i^-$ and $\lambda_j$ are unknowns.

Thus, the extent of input congestion related to input $x_{i0}$ can be described as

$$s_i^c = s_i^{-*} - \delta_i^{-*}, i=1,...,m \quad (4)$$

It should be noted that the importance of each input component is identical in the CTT model, so the units of inputs should be selected carefully in real applications.

**The third approach** is the WY-TS model proposed by (Wei and Yan, 2004; Tone and Sahoo, 2004). In this approach, we first calculate the PTE via the following Model (5):

$$\max_{\lambda_j, \pi} \ \pi$$

$$s.t. \begin{cases} \sum_{j=1}^{n} \lambda_j x_{ij} = x_{i0}, i = 1, \ldots, m \\ \sum_{j=1}^{n} \lambda_j y_{rj} \geq \pi y_{r0}, r = 1, \ldots, s \\ \sum_{j=1}^{n} \lambda_j = 1 \\ \lambda_j \geq 0, j = 1, \ldots, n, \pi \ free \end{cases} \tag{5}$$

The PPS implied in Model (5) is as follows.

$$P_{convex}(X,Y) = \left\{ (X,Y) \Big| \sum_{j=1}^{n} \lambda_j X_j = X, \sum_{j=1}^{n} \lambda_j Y_j \geq Y, \sum_{j=1}^{n} \lambda_j = 1, \lambda_j \geq 0, j = 1, \ldots, n \right\}$$

where $X_j = (x_{1j}, x_{2j}, \ldots, x_{mj})$ and $Y_j = (y_{1j}, y_{2j}, \ldots, y_{sj})$ denote the inputs and outputs, respectively. We can observe that this PPS allows the frontier to bend downward, on which the scale elasticity (SE) of the DMU is negative, i.e., congestion effect occurs. Let $\pi^*$ represent the optimal objective value of Model (5). The ratio $\varphi = \pi^* / \theta^*$ can be used as the measurement of the extent of congestion effect.

The dual model of Model (5) is exhibited as Model (6):

$$\max_{u_r, v_i, v_0} \ \sum_{r=1}^{s} u_r y_{r0}$$

$$s.t. \begin{cases} \sum_{r=1}^{s} u_r y_{rj} - \sum_{i=1}^{m} v_i x_{ij} - v_0 \leq 0, j = 1, \ldots, n \\ \sum_{i=1}^{m} v_i x_{i0} + v_0 = 1 \\ u_r \geq 0, r = 1, \ldots, s; v_i \ free, \ i = 0, 1, \ldots, m \end{cases} \tag{6}$$

Model (6) is similar to the classic BCC model proposed by Banker et al. (1984) except that the constraints $v_i, i = 1, \ldots, m$ are not non-negative but free. This fact indicates that the frontier of PPS may bend downward, i.e., the congestion effect occurs.

In particular, Tone and Sahoo (2004) proposed the definitions on strong congestion and weak congestion as follows.

**Definition 2 (Strong Congestion):** A DMU $(X_0, Y_0)$ is strongly congested if there exists an activity $(\tilde{X}_0, \tilde{Y}_0) \in P_{convex}(X,Y)$ such that $\tilde{X}_0 = \alpha X_0$ ( $0 < \alpha < 1$ ) and $\tilde{Y}_0 \geq \beta Y_0$ (with $\beta > 1$ ), where $(X_0, Y_0) \in P_{convex}(X,Y)$ and is strongly efficient.

**Definition 3 (Weak Congestion):** A DMU is (weakly) congested if it is strongly efficient with respect to $P_{convex}(X,Y)$ and there exists an activity in $P_{convex}(X,Y)$ that uses less resources in one or more inputs for making more products in one or more outputs.

Tone and Sahoo (2004) proposed the method to determine the DMU's status of strong and weak congestion.

**Step 1:** Solve the output-based BCC model (Model (1)), and then we have:

(a) If $\theta^* = 1, S^{-*} = 0, S^{+*} = 0$, then DMU $(X_0, Y_0)$ is BCC-efficient and not congested.

(b) If $\theta^* = 1, S^{-*} \neq 0, S^{+*} = 0$, then DMU $(X_0, Y_0)$ is BCC-inefficient.

(c) If $(\theta^* = 1, S^{+*} \neq 0)$ or $\theta^* > 1$, then DMU $(X_0, Y_0)$ displays congestion. Go to **Step 2**.

**Step 2:** Calculate the following Model (7):

$$\bar{\omega} = \max_{u_r, v_i, \mu_0} \mu_0$$

$$s.t. \begin{cases} \sum_{r=1}^s u_r y_{rj} - \sum_{i=1}^m v_i x_{ij} + \mu_0 \leq 0, j = 1,...,n \\ \sum_{r=1}^s u_r y_{r0} - \sum_{i=1}^m v_i x_{i0} + \mu_0 = 0 \\ \sum_{r=1}^s u_r y_{r0} = 1 \\ u_r \geq 0, r = 1,...,s; v_i \, free, i = 1,...,m; \mu_0 \, free \end{cases} \quad (7)$$

We let $\bar{\rho} = 1 + \bar{\omega}$. If $\bar{\rho} < 0$, then DMU $(X_0, Y_0)$ is strongly congested. If $\bar{\rho} \geq 0$, then DMU $(X_0, Y_0)$ is weakly congested.

Kao (2010) noted that the congestion measurement in Färe et al. (1985) is able to detect congestion for cases in which the FGL methods fails, and it is superior in differentiating the PTE with the congestion effect. Its deficiency is that it is unable to identify the excessive amount of each input that causes congestion.

## 3 Definitions of Directional Congestion

In this section, we use the same PPS $P_{convex}(X,Y)$ implied in the WY-TS model. Based on $P_{convex}(X,Y)$, we give the following definitions.

**Definition 4 (Weakly and Strongly efficient frontiers):** The weakly and strongly efficient frontiers of $P_{convex}(X,Y)$ can be defined as follows.

(1) Weakly efficient frontier:

$$EF_{convex\_weak} = \left\{(X,Y) \in P_{convex}(X,Y) \middle| \text{there is no } (\overline{X},\overline{Y}) \in P_{convex}(X,Y) \text{ such that } (-\overline{X},\overline{Y}) > (-X,Y) \right\} \quad (8a)$$

(2) Strongly efficient frontier:

$$EF_{convex\_strong} = \left\{(X,Y) \in P_{convex}(X,Y) \middle| \text{there is no } (\overline{X},\overline{Y}) \in P_{convex}(X,Y) \text{ such that } (-\overline{X},\overline{Y}) \geq (-X,Y) \text{ and } (\overline{X},\overline{Y}) \neq (X,Y) \right\} (8b)$$

**Definition 5 (Supporting Hyperplane):** A supporting hyperplane of $P_{convex}(X,Y)$ can be defined as follows. If a hyperplane

$$H(V,U,u_0) = \left\{(X,Y) \mid U^T Y - V^T X + u_0 = 0\right\} \quad (9)$$

satisfies (1) $(X_0,Y_0) \in \left\{(X,Y) \mid U^T Y - V^T X + u_0 = 0\right\} \bigcap P_{convex}(X,Y)$;

(2) For all $(X,Y) \in P_{convex}(X,Y)$, we have $U^T Y - V^T X + u_0 \leq 0$;

(3) $(V,U) \neq (\mathbf{0},\mathbf{0})$

Then, we say that $H(V,U,\mu_0)$ is a supporting hyperplane of $P_{convex}(X,Y)$ on the DMU $(X_0,Y_0)$, which is referred to as $H(V,U,\mu_0)\big|_{(X_0,Y_0)}$.

**Definition 6 (Face):** A subset of $P_{convex}(X,Y)$ is referred to as a "Face" on the DMU $(X_0,Y_0)$ if there exists a supporting hyperplane $H(V,U,\mu_0)\big|_{(X_0,Y_0)}$ such that the subset is identical to an intersection between $P_{convex}(X,Y)$ and $H(V,U,\mu_0)\big|_{(X_0,Y_0)}$.

We let $(\omega_1,...,\omega_m)^T$ ( $\omega_i \geq 0, i = 1,...,m$ ) and $(\delta_1,...,\delta_s)^T$ ( $\delta_r \geq 0, r = 1,...,s$ ) represent the input and output directions, respectively, and satisfy $\sum_{i=1}^m \omega_i = m; \sum_{r=1}^s \delta_r = s$ . We let $\Omega_0 = diag\{1 + \omega_1 t_0,...,1 + \omega_m t_0\}$ and $\Phi_0 = diag\{1 + \delta_1 \beta_0,...,1 + \delta_s \beta_0\}$, where $diag\{\cdot\}$ denotes the diagonal matrix, the constants $t_0$ and $\beta_0$ are the input and output scaling factors, respectively, and satisfy $t_0 \times \beta_0 < 0$ . Then, we give the following definitions on the directional congestion and local directional congestion.

**Definition 7 (Directional Congestion):** A DMU $(X_0,Y_0) \in P_{convex}(X,Y)$ is directionally congested if it is strongly efficient with respect to $P_{convex}(X,Y)$ and there exists an activity $\left(\tilde{X},\tilde{Y}\right)=\left(\Omega_0 X_0,\Phi_0 Y_0\right)\in P_{convex}\left(X,Y\right)$ that is also strongly efficient and uses less (more) resources in inputs in the direction $\left(\omega_1,...,\omega_m\right)^T$ for making more (less) products in outputs in the direction $\left(\delta_1,...,\delta_s\right)^T$.

**Definition 8 (Local Directional Congestion):** A DMU $\left(X_0,Y_0\right)\in P_{convex}\left(X,Y\right)$ is directionally congested if it is strongly efficient with respect to $P_{convex}\left(X,Y\right)$ and there exists an activity $\left(\tilde{X},\tilde{Y}\right)=\left(\Omega_0 X_0,\Phi_0 Y_0\right)\in P_{convex}\left(X,Y\right)$ on the same "Face" with DMU $\left(X_0,Y_0\right)$ that is also strongly efficient and uses less (more) resources in inputs in the direction $\left(\omega_1,...,\omega_m\right)^T$ for making more (less) products in outputs in the direction $\left(\delta_1,...,\delta_s\right)^T$.

**Theorem 1:** A DMU $\left(X_0,Y_0\right)\in P_{convex}\left(X,Y\right)$ is directionally congested if and only if it is locally directionally congested.

**Proof.** Please see Appendix A.

From **Theorem 1**, we can see that **Definition 7** is equivalent to **Definition 8**. Therefore, in this paper the directional congestion we discuss hereafter refers to the local directional congestion unless it is expressly stated.

From the definitions of directional congestion, we know that there are two possibilities for directional congestion: If a strongly efficient DMU $\left(X_0,Y_0\right)\in P_{convex}\left(X,Y\right)$ is directionally congested, there exists at least an activity $\left(\tilde{X},\tilde{Y}\right)=\left(\Omega_0 X_0,\Phi_0 Y_0\right)\in P_{convex}\left(X,Y\right)$ that is also strongly efficient and (1) either uses less resources in inputs in the direction $\left(\omega_1,...,\omega_m\right)^T$ for making more products in outputs in the direction $\left(\delta_1,...,\delta_s\right)^T$; or (2) uses more resources in inputs in the direction $\left(\omega_1,...,\omega_m\right)^T$ for making less products in outputs in the direction $\left(\delta_1,...,\delta_s\right)^T$. Take a production process with a single input and single output as an example. See Figure 2 for details.

In Figure 2, Point E is on the strongly efficient frontier $EF_{convex\_strong}$ with respect to $P_{convex}\left(X,Y\right)$. We can observe that there exists at least an activity that is strongly efficient and uses fewer inputs making more outputs at the left side of Point E. At the same time, there also exists at least an activity that is strongly efficient and uses more inputs making fewer outputs at the right side of Point E.

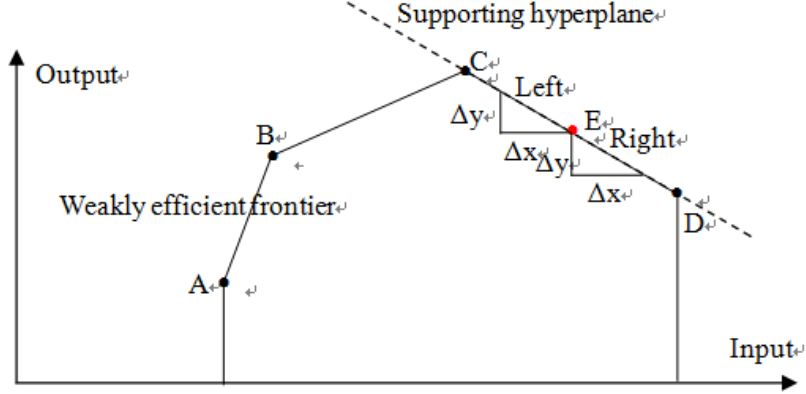

**Figure 2: Congestion effect on Point E**[2]

There may exist some strongly efficient $\left(X_0, Y_0\right) \in P_{convex}\left(X, Y\right)$ whose input(s) cannot be further reduced in a direction of $\left(\omega_1, \omega_2, ..., \omega_m\right)^T$ regardless of the output direction $\left(\delta_1, \delta_2, ..., \delta_s\right)^T$. These strongly efficient DMUs are defined as the directional smallest scale size (DSSS; see **Definition 9**). For example, see Point A in Figure 2. Analogously, there also may exist some strongly efficient $\left(X_0, Y_0\right) \in P_{convex}\left(X, Y\right)$ whose input(s) cannot be further expanded in a direction of $\left(\omega_1, \omega_2, ..., \omega_m\right)^T$ regardless of the output direction $\left(\delta_1, \delta_2, ..., \delta_s\right)^T$. These strongly efficient DMUs are defined as the directional largest scale size (DLSS; see **Definition 10**). For example, see Point D in **Figure 2**.

We let $\Omega_t = diag\left\{1 + \omega_1 t, ..., 1 + \omega_m t\right\}$ and $\Phi_\beta = diag\left\{1 + \delta_1 \beta, ..., 1 + \delta_s \beta\right\}$, where the variables $t$ and $\beta$ are input and output scaling factors, respectively.

**Definition 9 (Directional Smallest Scale Size, DSSS):** The strongly efficient $\left(X_0, Y_0\right) \in P_{convex}\left(X, Y\right)$ is of the directional smallest scale size if and only if $\left(\Omega_t X_0, \Phi_\beta Y_0\right) \notin P_{convex}\left(X, Y\right)$ for any $\beta$ and $t < 0$.

**Definition 10 (Directional Largest Scale Size, DLSS):** The strongly efficient $\left(X_0, Y_0\right) \in P_{convex}\left(X, Y\right)$ is of the directional largest scale size if and only if $\left(\Omega_t X_0, \Phi_\beta Y_0\right) \notin P_{convex}\left(X, Y\right)$ for any $\beta$ and $t > 0$.

Please note that the definitions of the directional smallest and largest scale size in our paper are compatible with the definition of extreme scale size in Banker and Thrall (1992) in the sense of the radial direction.

---

[2] Please note: In the case of a single input and a single output, directional congestion is equivalent to classic congestion.

# 4 Measurement of Directional Congestion

The DMUs are often divided into two categories in the measurement of congestion using the DEA models. The two categories are referred to as the strongly efficient DMUs[3] on the efficient frontier and the weakly efficient or inefficient DMUs. The congestion of weakly efficient or inefficient DMUs can be measured through their projections onto the strongly efficient frontier. This paper follows the above two categories and conducts the directional congestion measurement based on the $P_{convex}(X,Y)$, focusing on the directional congestion measurement of the strongly efficient DMUs. Similar to Yang et al. (2014), in this paper we will use the **Finite Difference Method (FDM)** and the **Upper and Lower Bounds Method (ULBM)** to measure the directional congestion effect.

## 4.1 Finite Difference Method (FDM)

Golany and Yu (1997) used the FDM to estimate the RTS for each DMU by testing the existence of solutions in four regions defined in the neighborhood of the analyzed unit. They provided a procedure to determine the RTS to the "right" and "left" of the DMU being evaluated. Rosen et al. (1998) estimated the directional derivative of DMUs on strongly efficient frontiers using the FDM. In this subsection, we will also use the FDM to determine the directional congestion effect of strongly efficient DMUs on the efficient frontier of $P_{convex}(X,Y)$.

### 4.1.1 Directional congestion of strongly efficient DMUs

It is well-known that the weakly or strongly efficient frontier of the BCC-DEA is piecewise linear. Thus, we determine the directional congestion to the "right" and "left" of the DMU being evaluated. Figure 2 exhibits the directional congestion to the "right" and "left" of the point E, which is on the strongly efficient frontier $EF_{convex\_strong}$.

#### 4.1.1.1 Directional congestion to the "right" of strongly efficient DMUs

First, we need to determine whether the strongly efficient DMU $(X_0,Y_0)$ is of the directional smallest scale

---

size. According to **Definition 10**, we consider the following Model (10):

$$\max_{\beta, \lambda_j, \eta} \eta$$

$$s.t. \begin{cases} \sum_{j=1}^{n} \lambda_j x_{ij} = \left(1 + \omega_i \eta\right) x_{i0}, i = 1,...,m \\ \sum_{j=1}^{n} \lambda_j y_{rj} \geq \left(1 + \delta_r \beta\right) y_{r0}, r = 1,...,s \\ \sum_{j=1}^{n} \lambda_j = 1, \lambda_j \geq 0, j = 1,...,n \\ \eta \geq 0; \beta \text{ free} \end{cases} \quad (10)$$

**Theorem 2.** The optimal objective value $\eta^*$ of Model (10) is zero if and only if the strongly efficient DMU $\left(X_0, Y_0\right)$ is of the directional largest scale size.

**Proof.** Please see Appendix A.

We first discuss the case in which the strongly efficient DMU $\left(X_0, Y_0\right)$ is not of the directional largest scale size in the direction of $\left(\omega_1, \omega_2,...,\omega_m\right)^T$ and $\left(\delta_1, \delta_2,...,\delta_s\right)^T$.

Based on the **Definition 7** and the FDM proposed by (Rosen et al., 1998; Golany and Yu, 1997), let $t_{right} > 0$. We have the following Model (11) to determine the right-hand directional congestion effect:

$$\max_{\beta, \lambda_j} \xi = \beta / t_{right}$$

$$s.t. \begin{cases} \sum_{j=1}^{n} \lambda_j x_{ij} = \left(1 + \omega_i t_{right}\right) x_{i0}, i = 1,...,m \\ \sum_{j=1}^{n} \lambda_j y_{rj} \geq \left(1 + \delta_r \beta\right) y_{r0}, r = 1,...,s \\ \sum_{j=1}^{n} \lambda_j = 1, \lambda_j \geq 0, j = 1,...,n \end{cases} \quad (11)$$

where $\delta_r \geq 0, r = 1,...,s$ and $\omega_i \geq 0, i = 1,...,m$ represent the direction factors of the inputs and outputs, respectively, and satisfy $\sum_{r=1}^{s} \delta_r = s; \sum_{i=1}^{m} \omega_i = m$. Model (11) appears to be a nonlinear programming. However, as observed below, its objective is independent of $t_{right} > 0$ after $t_{right}$ is sufficiently small. Thus, we will understand that $t_{right}$ is a small positive quantity, which represents the amount of the directional change of inputs. Variable $\beta$ represents the amount of the directional change of outputs. Then, it actually becomes a linear programming.

Let $\Omega_t^+ = diag\left\{1 + \omega_1 t_{right},...,1 + \omega_m t_{right}\right\}$ and $\Phi_{\beta^*}^+ = diag\left\{1 + \delta_1 \beta^*,...,1 + \delta_s \beta^*\right\}$, where $\beta^*$ is the optimal solution of Model (11). We have the following **Theorem 3**:

**Theorem 3.** There exists $t_0^+ > 0$ satisfying (1) when $t_{right} \in \left(0, t_0^+\right]$ and $\left(\Omega_t^+ X_0, \Phi_{\beta^*}^+ Y_0\right) \in P_{convex}\left(X, Y\right)$, $\left(\Omega_t^+ X_0, \Phi_{\beta^*}^+ Y_0\right)$ is located on the weakly efficient frontier $EF_{convex\_weak}$ and (2) when $t_{right} \in \left(0, t_0^+\right]$, $\left(X_0, Y_0\right)$ and $\left(\Omega_t^+ X_0, \Phi_{\beta^*}^+ Y_0\right)$ have the same supporting hyperplane, which may be different for different $t_{right}$.

**Proof.** Please see Appendix A.

**Definition 11:** Let $\xi^*\left(X_0, Y_0\right)$ be the optimal objective of Model (11). Thus, we can determine the directional congestion to the "right" of DMU $\left(X_0, Y_0\right)$ as follows: If $\xi^*\left(X_0, Y_0\right) < 0$ holds, then DMU $\left(X_0, Y_0\right)$ is directionally congested in the direction of $\left(\omega_1, \omega_2, ..., \omega_m\right)^T$ and $\left(\delta_1, \delta_2, ..., \delta_s\right)^T$.

Here, we demonstrate that for very small positive $t_{right}$, the objective value in Model (11) is a constant for any given input and output direction.

**Theorem 4.** There exists a sufficiently small quantity $t_0^+ > 0$ such that the optimal value of Model (11) is constant for all $t_{right} \in \left(0, t_0^+\right]$.

**Proof.** Please see Appendix A.

Now, we discuss how to select $t_{right}$ in practice. Consider the following Model (12):

$$\max_{\beta, \lambda_j} \beta$$
$$s.t. \begin{cases} \sum_{j=1}^n \lambda_j x_{ij} = \left(1 + \omega_i t_{right}\right) x_{i0}, i = 1, ..., m \\ \sum_{j=1}^n \lambda_j y_{rj} \geq \left(1 + \delta_r \beta\right) y_{r0}, r = 1, ..., s \\ \sum_{j=1}^n \lambda_j = 1, \lambda_j \geq 0, j = 1, ..., n \end{cases} \quad (12)$$

Let $\left(\lambda_j^*, \beta^*\right)$ be optimal solutions of Model (12). Consider the following Model (13):

$$\max_{U, V, \mu_0} \varphi = U^T Y_0 + \mu_0$$
$$s.t. \begin{cases} U^T Y_j - V^T X_j + \mu_0 \leq 0, j = 1, ..., n \\ V^T X_0 = 1 \\ U^T \Phi_{\beta^*}^+ Y_0 - V^T \Omega_t^+ X_0 + \mu_0 = 0 \\ U \geq \mathbf{0}; V, \mu_0 \text{ free} \end{cases} \quad (13)$$

From the proof of **Theorem 4**, if the optimal objective value of Model (13) satisfies $\varphi^* = 1$, the positive constant $t_{right}$ is sufficiently small to guarantee that both $(X_0, Y_0)$ and $(\Omega_t^+ X_0, \Phi_{\beta'}^+ Y_0)$ are located on the weakly efficient frontier $EF_{convex\_weak}$, and they have the same supporting hyperplane. Thus, in practice, we first select a small number $t_{right}$ in (12) and solve (13) to determine if the optimal is the unit. If not, we will attempt smaller numbers. From the continuity, it will be the unit when $t_{right}$ is sufficiently small.

Now, we turn to the case in which the strongly efficient DMU $(X_0, Y_0)$ is of the directional largest scale size in the direction of $(\omega_1, \omega_2, ..., \omega_m)^T$ and $(\delta_1, \delta_2, ..., \delta_s)^T$. In this case, we cannot find a feasible solution in Model (11) when $t_{right} > 0$ is a small positive constant. Thus, we provide the following **Definition 12** to address the right-hand directional RTS of $(X_0, Y_0)$:

**Definition 12:** If the strongly efficient $(X_0, Y_0)$ is of the directional largest scale size in the direction of $(\omega_1, \omega_2, ..., \omega_m)^T$ and $(\delta_1, \delta_2, ..., \delta_s)^T$, there is no data to determine the congestion effect on the right-hand side of DMU $(X_0, Y_0)$.

#### 4.1.1.2 Directional congestion to the "left" of strongly efficient DMUs

First, we need to determine whether the strongly efficient DMU $(X_0, Y_0)$ is of the directional smallest scale size. According to **Definition 9**, we consider the following Model (14):

$$\max_{\beta, \lambda, \eta} \eta$$

$$s.t. \begin{cases} \sum_{j=1}^n \lambda_j x_{ij} = (1 - \omega_i \eta) x_{i0}, i = 1, ..., m \\ \sum_{j=1}^n \lambda_j y_{rj} \geq (1 - \delta_r \beta) y_{r0}, r = 1, ..., s \\ \sum_{j=1}^n \lambda_j = 1, \lambda_j \geq 0, j = 1, ..., n \\ \eta \geq 0; \beta \text{ free} \end{cases} \quad (14)$$

**Theorem 5.** The optimal objective value $\eta^*$ of Model (14) is zero if and only if the strongly efficient DMU $(X_0, Y_0)$ is of the directional smallest scale size.

**Proof.** Please see Appendix A.

We first discuss the case in which the strongly efficient DMU $(X_0, Y_0)$ is not of the directional smallest scale size in the direction of $(\omega_1, \omega_2, ..., \omega_m)^T$ and $(\delta_1, \delta_2, ..., \delta_s)^T$. Based on the **Definition 7** and the FDM proposed by (Rosen et al., 1998; Golany and Yu, 1997), we let $t_{left}$ be a small positive constant and have the following Model (15) to determine the left-hand directional congestion effect:

$$\min_{\beta,\lambda_j} \ \psi = \beta / t_{left}$$

$$s.t. \begin{cases} \sum_{j=1}^{n} \lambda_j x_{ij} = \left(1 - \omega_i t_{left}\right) x_{i0}, i = 1,...,m \\ \sum_{j=1}^{n} \lambda_j y_{rj} \geq \left(1 - \delta_r \beta\right) y_{r0}, r = 1,...,s \\ \sum_{j=1}^{n} \lambda_j = 1, \lambda_j \geq 0, j = 1,...,n \end{cases} \qquad (15)$$

where $\delta_r \geq 0, r = 1,...,s$ and $\omega_i \geq 0, i = 1,...,m$ represent the direction factors of the inputs and outputs, respectively and satisfy $\sum_{r=1}^{s} \delta_r = s; \sum_{i=1}^{m} \omega_i = m$. Constant $t_{left}$ is a small positive quantity that represents the amount of the directional change of inputs. Variable $\beta$ represents the amount of the directional change of outputs.

We let $\Omega_t^- = diag\left\{1 - \omega_1 t_{left},...,1 - \omega_m t_{left}\right\}$ and $\Phi_{\beta^*}^- = diag\left\{1 - \delta_1 \beta^*,...,1 - \delta_s \beta^*\right\}$, and $\beta^*$ is the optimal solution of Model (15). Thus, we have the following **Theorem 6**:

**Theorem 6.** There exists $t_0^- > 0$ satisfying (1) when $t_{left} \in \left(0, t_0^-\right]$ and $\left(\Omega_t^- X_0, \Phi_{\beta^*}^- Y_0\right) \in P_{convex}\left(X, Y\right)$, $\left(\Omega_t^- X_0, \Phi_{\beta^*}^- Y_0\right)$ is located on the weakly efficient frontier $EF_{convex\_weak}$ and (2) when $t_{left} \in \left(0, t_0^-\right]$, $\left(X_0, Y_0\right)$ and $\left(\Omega_t^- X_0, \Phi_{\beta^*}^- Y_0\right)$ have the same supporting hyperplane, which may be different for different $t_{left}$.

**Proof.** Please see Appendix A.

**Definition 13:** We let $\psi^*\left(X_0, Y_0\right)$ be the optimal objective of Model (15). Accordingly, we can determine the directional congestion to the "left" of the DMU $\left(X_0, Y_0\right)$ as follows: if $\psi^*\left(X_0, Y_0\right) < 0$ holds, then $DMU\left(X_0, Y_0\right)$ is directionally congested in the direction of $\left(\omega_1, \omega_2,...,\omega_m\right)^T$ and $\left(\delta_1, \delta_2,...,\delta_s\right)^T$.

Next, we discuss how to choose $t_{left}$. Again, we first have

**Theorem 7.** There exists a sufficiently small quantity $t_0^- > 0$ such that the optimal value of Model (15) is constant for all $t_{left} \in \left(0, t_0^-\right]$.

**Proof.** Please see Appendix A.

Now, again, consider the following Model (16):

$$\max_{\beta,\lambda_j} \beta$$

$$s.t. \begin{cases} \sum_{j=1}^{n} \lambda_j x_{ij} = \left(1 - \omega_i t_{left}\right) x_{i0}, i = 1,...,m \\ \sum_{j=1}^{n} \lambda_j y_{rj} \geq \left(1 - \delta_r \beta\right) y_{r0}, r = 1,...,s \\ \sum_{j=1}^{n} \lambda_j = 1, \lambda_j \geq 0, j = 1,...,n \end{cases} \quad (16)$$

Let $\left(\lambda_j^*, \beta^*\right)$ be the optimal solutions of Model (16) and $\Omega_t^- = diag\left\{1 - \omega_1 t_{left},...,1 - \omega_m t_{left}\right\}$ and $\Phi_{\beta^*}^- = diag\left\{1 - \delta_1 \beta^*,...,1 - \delta_s \beta^*\right\}$. Consider the following Model (17):

$$\max_{U,V,\mu_0} \varphi = U^T Y_0 + \mu_0$$

$$s.t. \begin{cases} U^T Y_j - V^T X_j + \mu_0 \leq 0, j = 1,...,n \\ V^T X_0 = 1 \\ U^T \Phi_t^- Y_0 - V^T \Omega_{\beta^*}^- X_0 + \mu_0 = 0 \\ U \geq \mathbf{0}, V, \mu_0 \text{ free} \end{cases} \quad (17)$$

Again, if the optimal objective value of Model (17) satisfies $\varphi^* = 1$, $t_{left} > 0$ is a sufficiently small constant to guarantee that both $\left(X_0, Y_0\right)$ and $\left(\Omega_t^- X_0, \Phi_{\beta^*}^- Y_0\right)$ are located on the weakly efficient frontier $EF_{convex\_weak}$, and they have the same supporting hyperplane. Thus, we will select $t_{left} > 0$ similarly.

Now, we turn to the case in which the strongly efficient $\left(Y_0, X_0\right)$ is of the directional smallest scale size in the direction of $\left(\omega_1, \omega_2,...,\omega_m\right)^T$ and $\left(\delta_1, \delta_2,...,\delta_s\right)^T$. In this case, we cannot find a feasible solution in Model (15) when $t_{left} > 0$ is a small positive constant. Thus, we provide the following **Definition 14** to address the left-hand directional congestion of $\left(X_0, Y_0\right)$:

**Definition 14:** If strongly efficient $\left(X_0, Y_0\right)$ is of the directional smallest scale size in the direction of $\left(\omega_1, \omega_2,...,\omega_m\right)^T$ and $\left(\delta_1, \delta_2,...,\delta_s\right)^T$, then there is no data to determine the congestion effect at the left-hand side of $\left(X_0, Y_0\right)$.

### 4.1.1.3 A procedure for estimating directional congestion of strongly efficient DMUs

Based on the above analysis, we now propose a procedure for estimating the directional congestion of a strongly efficient DMU $\left(X_0, Y_0\right)$ as follows.

**Procedure 1:** The directional congestion to the "right" and "left" of a strongly efficient DMU $\left(X_0, Y_0\right)$ on the strongly efficient frontier $EF_{convex\_strong}$ with respect to $P_{convex}(X,Y)$ in the direction of $(\omega_1,\omega_2,...,\omega_m)^T$ and $(\delta_1,\delta_2,...,\delta_s)^T$ can be determined by:

**(a) Right-hand directional congestion**

**Step a-0:** Solve Model (10) to determine whether its optimal objective value is zero. If so, DMU$(X_0,Y_0)$ is of the directional largest scale size, and there is no data to determine the right-hand directional congestion. Otherwise, we have the following two steps:

**Step a-1:** Choose a sufficiently small quantity $t_{right}>0$, based on Model (12) and Model (13), to guarantee that both $(X_0,Y_0)$ and $\left(\Omega_t^+ X_0, \Phi_\beta^+ Y_0\right)$ are located on the weakly efficient frontier $EF_{convex\_weak}$, and they have the same supporting hyperplane.

**Step a-2:** Solve Model (11) to determine the directional congestion to the "right" of DMU$(X_0,Y_0)$ in the direction of $(\omega_1,\omega_2,...,\omega_m)^T$ and $(\delta_1,\delta_2,...,\delta_s)^T$: If $\xi^*(X_0,Y_0)<0$, then DMU$(X_0,Y_0)$ is directionally congested.

**(b) Left-hand directional congestion**

**Step b-0:** Solve Model (14) to determine whether its optimal objective value is zero. If so, DMU$(X_0,Y_0)$ is of the directional smallest scale size, and there is no data to determine the left-hand directional congestion. Otherwise, we have the following two steps:

**Step b-1:** Choose a sufficiently small quantity $t_{left}>0$, based on Model (16) and Model (17), to guarantee that both $(X_0,Y_0)$ and $\left(\Omega_t^- X_0, \Phi_\beta^- Y_0\right)$ are located on the weakly efficient frontier $EF_{convex\_weak}$ and that they have the same supporting hyperplane.

**Step b-2:** Solve Model (15) to determine the directional congestion to the "right" of DMU$(X_0,Y_0)$ in the direction of $(\omega_1,\omega_2,...,\omega_m)^T$ and $(\delta_1,\delta_2,...,\delta_s)^T$: If $\psi^*(X_0,Y_0)<0$, then DMU$(X_0,Y_0)$ is directionally congested.

**4.1.2 Directional congestion of inefficient or weakly efficient DMUs**

For estimating the directional congestion to the "right" and "left" of inefficient or weakly efficient DMUs, we can perform the following two steps:

**Step 1:** First, we project the inefficient or weakly efficient DMUs onto the strongly efficient frontier $EF_{convex\_strong}$ using the following Model (18):

$$\max_{\lambda_j, \theta_0, s_r^-, s_r^+} \theta_0$$

$$s.t. \begin{cases} \sum_j \lambda_j x_{ij} = x_{i0}, i = 1,...,m \\ \sum_j \lambda_j y_{rj} - s_r^+ = \theta_0 y_{r0}, r = 1,...,s \\ \sum_j \lambda_j = 1 \\ \lambda_j, j = 1,...,n; s_r^-, s_r^+ \geq 0, r = 1,...,s; i = 1,...,m \end{cases} \quad (18)$$

, and the formula for projection is the following Equation (19):

$$\begin{cases} \tilde{x}_{i0} \leftarrow x_{i0}, i = 1,...,m \\ \tilde{y}_{r0} \leftarrow \theta_0^* y_{r0} + s_r^{+*}, r = 1,...,s \end{cases} \quad (19)$$

where ($\theta^*, s_r^{-*}$) is the optimal solution of Model (18).

**Step 2:** When we determine the projected points on the strongly efficient frontier $EF_{convex\_strong}$, we can estimate the directional congestion to the "right" and "left" for the inefficient or weakly efficient DMUs using **Procedure 1** in Section 4.1.1.3.

**Remark 1.** *Please note that different projections may generate different results (see, e.g., Sueyoshi and Sekitani, 2009). In this paper, we focus mainly on the detection of the directional congestion of the DMUs on the strongly efficient frontier* $EF_{convex\_strong}$ *instead of the inefficient or weakly efficient ones.*

### 4.2 Upper and Lower Bounds Method (ULBM)

In this subsection, we mainly discuss the directional congestion of strongly efficient DMU $(X_0, Y_0)$ on the strongly efficient frontier $EF_{convex\_strong}$ with respect to $P_{convex}(X, Y)$ in the direction of $(\omega_1, \omega_2, ..., \omega_m)^T$ and $(\delta_1, \delta_2, ..., \delta_s)^T$ using another approach, i.e., the upper and lower bound (ULBM) method.

According to **Theorem 3** and **Theorem 6**, we know that the following formulas hold when $t_{right} > 0$ and $t_{left} > 0$ are sufficiently small positive constants in cases in which DMU $(X_0, Y_0)$ is neither the DLSZ nor DSSZ:

$$\frac{\beta^*}{t_{right}} = \frac{\sum_{i=1}^m v_i^* \omega_i x_{i0}}{\sum_{r=1}^s u_r^* \delta_r y_{r0}} \quad \text{and} \quad \frac{\beta^*}{t_{left}} = \frac{\sum_{i=1}^m v_i^* \omega_i x_{i0}}{\sum_{r=1}^s u_r^* \delta_r y_{r0}} \quad (20)$$

where $U^* = \left(u_1^*, u_2^*, ..., u_s^*\right)^T$ and $V^* = \left(v_1^*, v_2^*, ..., v_m^*\right)^T$ are the optimal solutions of Model (13) and Model (17), respectively.

Thus, we can use the following Model (21) to calculate the upper and lower bounds of the directional SE and, subsequently, determine the directional congestion of DMU $\left(X_0, Y_0\right)$.

$$
\bar{\rho}\left(\underline{\rho}\right) = \max_{u_r, v_i, \mu_0}\left(\min_{u_r, v_i, \mu_0}\right)\frac{\sum_{i=1}^m v_i \omega_i x_{i0}}{\sum_{r=1}^s u_r \delta_r y_{r0}}
$$

$$
s.t.\begin{cases} \sum_{r=1}^s u_r y_{rj} - \sum_{i=1}^m v_i x_{ij} + \mu_0 \le 0, j = 1,...,n \\ \sum_{r=1}^s u_r y_{r0} - \sum_{i=1}^m v_i x_{i0} + \mu_0 = 0 \\ \sum_{i=1}^m v_i x_{i0} - \mu_0 = 1 \\ u_r \ge 0, r = 1,...,s, v_i, \mu_0 \text{ free}, i = 1,...,m \end{cases} \quad (21)
$$

Therefore, in the case in which DMU $\left(X_0, Y_0\right)$ is neither the DLSS nor DSSS, we have the following theorem:

**Theorem 8.** Suppose that DMU $\left(X_0, Y_0\right)$ is neither the DLSS nor DSSS. The upper and lower bounds ($\bar{\rho}$ $\left(X_0, Y_0\right)$ and $\underline{\rho}\left(X_0, Y_0\right)$) in Model (21) are equal to the optimal objective value $\psi^*\left(X_0, Y_0\right)$ in Model (15) and the optimal objective value $\xi^*\left(X_0, Y_0\right)$ in Model (11).

**Proof.** Please see Appendix A.

**Theorem 9.** (1) If the maximal optimal objective value $\bar{\rho}\left(X_0, Y_0\right)$ of Model (21) is unbounded ($+\infty$), the strongly efficient DMU $\left(X_0, Y_0\right)$ is of the directional smallest scale size. (2) If the minimal optimal objective value $\underline{\rho}\left(X_0, Y_0\right)$ of Model (21) is unbounded ($-\infty$), the strongly efficient DMU $\left(X_0, Y_0\right)$ is of the directional largest scale size.

**Proof.** Please see Appendix A.

Therefore, we have the procedure similar to **Procedure 1** for estimating directional congestion.

Model (21) is a fractional programming that is difficult to solve, so we transform Model (21) into an equivalent linear programming through the Charnes-Cooper transformation (Charnes et al., 1962). First, we rewrite Model (21) using the following Model (22):

$$\bar{\rho}\left(\underline{\rho}\right) = \max_{U,V,\mu_0}\left(\min_{U,V,\mu_0}\right)\frac{V^T W X_0}{U^T \Delta Y_0}$$

$$s.t.\begin{cases} U^T Y_j - V^T X_j + \mu_0 \leq 0, j = 1,...,n \\ U^T Y_0 - V^T X_0 + \mu_0 = 0 \\ V^T X_0 - \mu_0 = 1 \\ U \geq \mathbf{0}, V, \mu_0 \text{ free} \end{cases} \tag{22}$$

where $U = \left(u_1, u_2,...,u_s\right)^T$ and $V = \left(v_1, v_2,...,v_m\right)^T$ are vectors of multipliers, and $\Delta = diag\left\{\delta_1, \delta_2,...,\delta_s\right\}$ and $W = diag\left\{\omega_1, \omega_2,...,\omega_m\right\}$ are matrixes of the directions of inputs and outputs.

We let

$$\begin{cases} \tau = 1/U^T \Delta Y_0 \\ \tau V^T W = \Gamma^T \\ \tau U^T \Delta = \Lambda^T \end{cases} \tag{23}$$

We assume the inverse matrix of $\Delta = diag\left\{\delta_1, \delta_2,...,\delta_s\right\}$ exists. Let $\omega_i \geq \varepsilon$, where $\varepsilon$ is a non-Archimedean construct to ensure the inverse matrix of $W$ exists. In this case, we have the following Equation (24) from Equation (23):

$$\begin{cases} \tau = 1/U^T \Delta Y_0 \\ V^T = \dfrac{1}{\tau}\Gamma^T W^{-1} \\ U^T = \dfrac{1}{\tau}\Lambda^T \Delta^{-1} \end{cases} \tag{24}$$

Thus, Model (22) can be translated into the following Model (25):

$$\bar{\rho}\left(\underline{\rho}\right) = \max_{\Gamma,\Lambda,\tau,\mu_0}\left(\min_{\Gamma,\Lambda,\tau,\mu_0}\right)\Gamma^T X_0$$

$$s.t.\begin{cases} \Lambda^T \Delta^{-1} Y_j - \Gamma^T W^{-1} X_j + \tau\mu_0 \leq 0, j = 1,...,n \\ \Lambda^T \Delta^{-1} Y_0 - \Gamma^T W^{-1} X_0 + \tau\mu_0 = 0 \\ \Gamma^T W^{-1} X_0 - \tau\mu_0 = \tau \\ \Lambda^T Y_0 = 1 \\ \Lambda \geq \mathbf{0}; \tau \geq 0; \Gamma, \mu_0 \text{ free} \end{cases} \tag{25}$$

We let $\tau\mu_0 = \mu_0'$, and Model (25) can be converted into the following linear programming:

$$\bar{\rho}\left(\underline{\rho}\right) = \max_{\Gamma, \Lambda, \tau, \mu_0'} \left( \min_{\Gamma, \Lambda, \tau, \mu_0'} \right) \Gamma^T X_0$$

$$s.t. \begin{cases} \Lambda^T \Delta^{-1} Y_j - \Gamma^T W^{-1} X_j + \mu_0' \leq 0, j = 1, \dots, n \\ \Lambda^T \Delta^{-1} Y_0 - \Gamma^T W^{-1} X_0 + \mu_0' = 0 \\ \Gamma^T W^{-1} X_0 - \mu_0' = \tau \\ \Lambda^T Y_0 = 1 \\ \Lambda \geq \mathbf{0}; \ \tau \geq 0; \ \Gamma, \mu_0' \text{ free} \end{cases} \quad (26)$$

Solving Model (26), we can obtain the optimal objective value of Model (21) or Model (22).

**4 .3 Some properties of directional congestion**

Using the definitions of strong and weak congestion in **Definition 2** and **Definition 3**, we can prove the following three theorems to address the relationships between directional congestion and strong (or weak) congestion.

**Theorem 10:** Strongly efficient DMU $\left( X_0, Y_0 \right) \in P_{convex} \left( X, Y \right)$ is left-hand directionally congested in the diagonal direction (i.e., $\omega_i = 1, i = 1, \dots, m$ ; $\delta_r = 1, r = 1, \dots, s$ ) if and only if it is strongly congested.

**Proof.** Please see Appendix A.

**Theorem 11:** Strongly efficient DMU $\left( X_0, Y_0 \right) \in P_{convex} \left( X, Y \right)$ is left-hand directionally congested in a certain direction if and only if it is weakly congested.

**Proof.** Please see Appendix A.

**Theorem 12:** If strongly efficient DMU $\left( X_0, Y_0 \right) \in P_{convex} \left( X, Y \right)$ is neither strongly congested nor weakly congested, there exists no directional congestion in any direction at its left-hand.

**Proof.** Please see Appendix A.

**Theorem 13:** Strongly efficient DMU $\left( X_0, Y_0 \right) \in P_{convex} \left( X, Y \right)$ displays congestion if it is directionally congested in a certain direction.

**Proof.** Please see Appendix A.

# 5 A Case Study

In this section, we conduct a case study to analyze the directional RTS of 16 basic research institutes in the Chinese Academy of Sciences (CAS) in 2010. Since the Pilot Project of Knowledge Innovation (KIPP) in 1998 at the CAS, institute evaluation has become increasingly important, and the requirements for the evaluation process have diversified. Since 2005, the CAS headquarters has built up the Comprehensive Quality Evaluation (CQE) system for institute evaluation in the CAS. The results of the evaluation are expressed as multi-dimensional feedback data and used as the tools to provide the basis of comprehensive analysis and decision-making and to provide institutes with targeted evaluation information and diagnostic comments. In the framework of the CQE, multiple inputs and outputs of the basic research institutes of the CAS are monitored using several quantitative indicators. In Liu et al. (2011), discussions are undertaken to select the suitable indicators for the DEA (data envelopment analysis)-based evaluations on basic research institutes in CAS. In this work, we will still use these indicators for analyzing the directional congestion for these basic research institutes. See Table 1 for details.

**Table 1: Input/output indicators**

| Indicators | Type | Units | Explanations |
|------------|------|-------|--------------|
| Staff | Input | Full Time Equivalent (FTE) | Full-time equivalent of full-time research staff |
| Res. Expen. | Input | RMB million | Amount of total income of each institute |
| SCI Pub. | Output | Number | Number of international papers indexed by the Web of Science from Thompson Reuters |
| High Pub. | Output | Number | Number of high-quality papers published in top research journals (e.g., journals with a top 15% impact factor) |
| Grad. Enroll. | Output | Number | Number of graduate student enrolment in 2009 |
| Exter. Fund | Output | RMB in million | Amount of external research funding from research contracts |

The data of these indicators (2 inputs and 4 outputs) are from the quantitative monitoring report in 2011 in the CAS and the Statistical Yearbook of CAS in 2011. See Table 2 for details.

First, we detect the congestion effect of 16 DMUs using the WY-TS model (Wei and Yan, 2004; Tone and Sahoo, 2004) based on the input-output data of these institutes by Model (18) and Equation (19). See Table B-1 in Appendix A for projections of these institutes. We can observe that the congestion effect occurs on $DMU_3$, $DMU_8$, $DMU_9$, $DMU_{10}$, $DMU_{11}$, $DMU_{12}$, $DMU_{15}$ and $DMU_{16}$. In the WY-TS model, we can further detect that these congested DMUs display weak or strong congestion. See Table 2 for details.

**Table 2. Input-output data and congestion effect using the WY-TS model.**

| DMUs | Outputs | | | | Inputs | | Congestion effect (WY-TS model) |
|------|---------|------|--------|-------------|-------|-----------|--------------------------------|
| | SCI | High | Grad.E | Exter.Fund. | Staff | Res.Expen. | |

| | Pub. | Pub. | nroll. | | | | PTE-Model (5)-$\eta^*$ | BCC-efficiency Model (1)-$\theta^*$ | Degree $\varphi=\eta^*/\theta^*$ | Strong or weak congestion |
|---|---|---|---|---|---|---|---|---|---|---|
| DMU$_1$ | 436 | 133 | 184 | 31.558 | 252 | 117.945 | 1 | 1 | 1 | No |
| DMU$_2$ | 243 | 127 | 43 | 15.3041 | 37 | 29.431 | 1 | 1 | 1 | No |
| DMU$_3$ | 164 | 70 | 89 | 33.8365 | 240 | 101.425 | 1.1835 | 1.4227 | 0.8319 | Weak |
| DMU$_4$ | 810 | 276 | 247 | 183.8434 | 356 | 368.483 | 1 | 1 | 1 | No |
| DMU$_5$ | 200 | 55 | 111 | 12.9342 | 310 | 195.862 | 1.9684 | 1.9684 | 1 | No |
| DMU$_6$ | 104 | 49 | 33 | 60.7366 | 201 | 188.829 | 1.6499 | 1.6499 | 1 | No |
| DMU$_7$ | 113 | 49 | 45 | 72.5368 | 157 | 131.301 | 1.0437 | 1.0437 | 1 | No |
| DMU$_8$ | 8 | 1 | 44 | 23.7015 | 236 | 77.439 | 1 | 1.9021 | 0.5257 | Weak |
| DMU$_9$ | 371 | 118 | 89 | 216.9885 | 805 | 396.905 | 1 | 1.2755 | 0.7840 | Strong |
| DMU$_{10}$ | 607 | 216 | 168 | 88.5561 | 886 | 411.539 | 1 | 1.5997 | 0.6251 | Strong |
| DMU$_{11}$ | 314 | 49 | 89 | 45.3597 | 623 | 221.428 | 1 | 2.1876 | 0.4571 | Weak |
| DMU$_{12}$ | 261 | 79 | 131 | 41.1156 | 560 | 264.341 | 1.4478 | 2.1500 | 0.6734 | Strong |
| DMU$_{13}$ | 627 | 168 | 346 | 645.4150 | 1344 | 900.509 | 1 | 1 | 1 | No |
| DMU$_{14}$ | 971 | 518 | 335 | 205.4528 | 508 | 344.312 | 1 | 1 | 1 | No |
| DMU$_{15}$ | 395 | 180 | 117 | 90.0373 | 380 | 161.331 | 1 | 1.1274 | 0.8870 | Weak |
| DMU$_{16}$ | 229 | 138 | 62 | 32.6111 | 132 | 83.972 | 1.3371 | 1.4111 | 0.9476 | Strong |

Data Source: (1) Monitoring data of the institutes in the CAS, 2011; (2) Statistical yearbook of CAS, 2011.
Note: These data were derived from these institutes in the period Jan.01, 2010~Dec.31, 2010.

Second, we can analyze the directional congestion effect of the above DMUs using our methods mentioned in Section 4.

**Step 1:** We use the output-based DEA model (Model (18)) to determine the strongly efficient frontier $EF_{convex\_strong}$ and the weakly efficient frontier $EF_{convex\_weak}$. Using Equation (19), we project the weakly efficient and inefficient DMUs onto the strongly efficient frontier $EF_{convex\_strong}$.

**Step 2:** According to Model (10) and Model (14), we can observe whether a DMU is DLSZ or DSSZ. Thus, we can use Model (11), Model (15) and Model (26) to determine the directional congestion on the right- and left-hand side of the DMUs. We take DMU$_1$ and DMU$_{15}$ as examples because they can exhibit certain interesting results. Without the loss of generality, we set the outputs direction as $\delta_1=\delta_2=\delta_3=\delta_4=1$, and we can compute the directional congestion effect of these two DMUs in different inputs directions using the FDM — we let $t_{right}=t_{left}=1E^{-6}$, which can pass the tests of Model (12)-Model (13) and Model (16)-Model (17) — and the ULBM, respectively. See Table 3 for detailed results.

**Table 3. Directional congestion effect of DMU$_1$ and DMU$_{15}$ in different inputs directions.**

| DMUs | $\omega_1$ | $\omega_2$ | $\xi^*$ (Right) | $\psi^*$ (Left) | $\underline{\rho}$ (Lower bound) | $\bar{\rho}$ (Upper bound) | Directional congestion effect (right) | Directional congestion effect (left) |
|------|-----------|-----------|-----------|-----------|-----------|-----------|-----------|-----------|
| | 0.3 | 1.7 | 0.07 | 4.64 | 0.07 | 4.64 | No | No |
| | 0.5 | 1.5 | 0.12 | 5.23 | 0.12 | 5.23 | No | No |
| | 0.7 | 1.3 | 0.17 | 7.32 | 0.17 | 7.32 | No | No |
| | 0.9 | 1.1 | 0.15 | 9.41 | 0.15 | 9.41 | No | No |
| DMU$_1$ | 1 | 1 | 0.14 | 10.46 | 0.14 | 10.46 | No | No |
| | 1.1 | 0.9 | 0.12 | 11.51 | 0.12 | 11.51 | No | No |
| | 1.3 | 0.7 | 0.10 | 13.60 | 0.10 | 13.60 | No | No |
| | 1.5 | 0.5 | 0.07 | 15.69 | 0.07 | 15.69 | No | No |
| | 1.7 | 0.3 | 0.04 | 17.78 | 0.04 | 17.78 | No | No |
| | 0.3 | 1.7 | 2.05 | 6.71 | 2.05 | 6.71 | No | No |
| | 0.5 | 1.5 | 1.72 | 5.03 | 1.72 | 5.03 | No | No |
| | 0.7 | 1.3 | 1.09 | 3.35 | 1.09 | 3.35 | No | No |
| | 0.9 | 1.1 | 0.37 | 1.67 | 0.37 | 1.67 | No | No |
| DMU$_{15}$ | 1 | 1 | 0 | 1.13 | 0 | 1.13 | No | No |
| | 1.1 | 0.9 | -0.55 | 0.85 | -0.55 | 0.85 | Yes | No |
| | 1.3 | 0.7 | -1.74 | 0.50 | -1.74 | 0.50 | Yes | No |
| | 1.5 | 0.5 | -3.40 | 0.16 | -3.40 | 0.16 | Yes | No |
| | 1.7 | 0.3 | -5.06 | -0.18 | -5.06 | -0.18 | Yes | Yes |

Based on the above analysis, we can find that the congestion effect occurs on DMU$_{15}$ when using the WY-TS model – in fact, this means that congestion can occur in the diagonal (strong congestion) or other (weak congestion) directions. From the directional congestion analysis, we know that whether the directional congestion effect occurs in certain directions (e.g., $\omega_1 = 1.7, \omega_2 = 0.3$; $\delta_1 = \delta_2 = \delta_3 = \delta_4 = 1$). On the other hand, for the same DMU, there are certain directions (e.g., $\omega_1 = 0.3, \omega_2 = 1.7$; $\delta_1 = \delta_2 = \delta_3 = \delta_4 = 1$) for which congestion does not occur. These are important details for institute planning. The congestion effect does not occur on DMU$_1$ in the WY-TS model, and the directional congestion effect does not occur to the left of DMU$_1$ as well. Similarly, we can analyze the directional congestion effect for other DMUs.

Furthermore, we can observe that there is weak congestion on DMU$_{15}$ according to the methods proposed by Tone and Sahoo (2004). The detailed procedure is as follows: We first solve the BCC model (Model (1)), and we have $\theta^*\left(DMU_{15}\right) = 1.1274 > 1$. Then, by solving Model (7), we have:

$$\bar{\rho} = 1 + \bar{\omega} = 1 + 0.1266 = 1.1266$$

Thus, we know DMU$_{15}$ is weakly congested. According to the model ([Congestion-weak]) in Tone and Sahoo (2004), we can easily have the direction of congestion as follows:

$$（\omega_1 = 2, \omega_2 = 0\,; \delta_1 = 1.9302, \delta_2 = 1.3981, \delta_3 = 0.6097, \delta_4 = 0.0621）$$

Therefore, we know that there exists directional congestion on DMU$_{15}$ at least in the above direction. In fact, from Table 3, we can observe that there is directional congestion in a certain direction, e.g., $（\omega_1 = 1.7, \omega_2 = 0.3\,;\ \ \delta_1 = \delta_2 = \delta_3 = \delta_4 = 1）$. In other words, if classic congestion occurs on a DMU, it is also directionally congested in a certain direction. This fact is consistent with the description in **Theorem 11**.

## 6 Conclusions

This paper proposes the definition of directional congestion and two methods named the FDM and ULBM are also presented to estimate the directional congestion to the right- and left-hand side of strongly efficient DMUs on $EF_{convex\_strong}$ with respect to $P_{convex}(X,Y)$. In addition, the relationships between the directional congestion and the classic strong (weak) congestion are explored in this paper. We find that: (1) The strongly efficient DMU is left-hand directionally congested in the diagonal direction if and only if it is strongly congested; (2) the strongly efficient DMU is left-hand directionally congested in a certain direction if and only if it is weakly congested; and (3) if the strongly efficient DMU is neither strongly congested nor weakly congested, no directional congestion in any direction on its left-hand side is observed. The directional congestion can be useful for decision-makers (DMs) to decide a rational combination of resources.

**Acknowledgments.** We would like to acknowledge the support of the National Natural Science Foundation of China (No. 71201158). The other supports of data and related materials from the Institute of Policy and Management, Chinese Academy of Sciences are also acknowledged.

# Appendix A

**Proof of Theorem 1.** Without loss of generality, we discuss the left-hand directional congestion and local directional congestion on strongly efficient DMU $(X_0, Y_0) \in P_{convex}(X, Y)$. See **Figure A-1**.

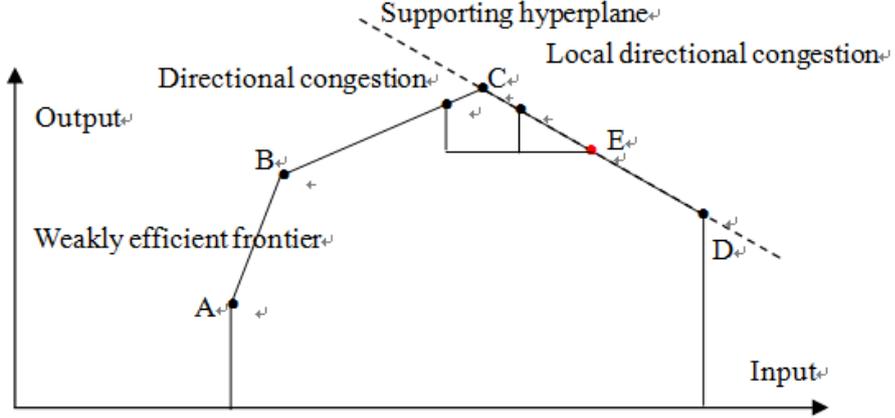

**Figure A-1: Left-hand directional congestion and local directional congestion on Point E.**

First, if directional congestion occurs on the left-hand side of strongly efficient DMU $(X_0, Y_0) \in P_{convex}(X, Y)$, there exists strongly efficient DMU $(\tilde{X}, \tilde{Y}) \in P_{convex}(X, Y)$ that satisfies $\tilde{X} \leq X_0$ and $\tilde{Y} \geq Y_0$. We let $F(X, Y) = 0$ represent the production function generated from the weak frontier with respect to $P_{convex}(X, Y)$. Because $P_{convex}(X, Y)$ is a convex set, we know that $F(X, Y) = 0$ is a piecewise linear concave function.

We let $\lambda \in [0, 1]$ and $\overline{X} = \lambda X_0 + (1 - \lambda)\tilde{X}$, where $\overline{Y}$ is the output produced by input $\overline{X}$ and satisfies $F(\overline{X}, \overline{Y}) = 0$. Because $F(X, Y) = 0$ is a piecewise linear concave function, we know that $\overline{X} = \lambda X_0 + (1 - \lambda)\tilde{X} \leq X_0$ and $\overline{Y} \geq \lambda Y_0 + (1 - \lambda)\tilde{Y} \geq Y_0$. Thus, we have $\overline{X} \leq X_0$ and $\overline{Y} \geq Y_0$. When $\lambda$ approaches 1, $(\overline{X}, \overline{Y})$ approaches $(X_0, Y_0)$. Thus, we can find $\lambda^*$ such that $(\overline{X}^*, \overline{Y}^*)$ and $(X_0, Y_0)$ are on the same "Face" of the weakly efficient frontier, where $\overline{X}^* = \lambda^* X_0 + (1 - \lambda^*)\tilde{X}$ and $\overline{Y}^*$ is the output produced by input $\overline{X}^*$ and satisfies $F(\overline{X}^*, \overline{Y}^*) = 0$. Thus, the local directional congestion occurs on the left-hand side of DMU $(X_0, Y_0)$.

Second, if the local directional congestion occurs on the left-hand side of DMU $(X_0, Y_0)$, we can easily observe that directional congestion occurs. **Q.E.D.**

**Proof of Theorem 2.** According to **Definition 10**, we can easily observe that **Theorem 2** holds. **Q.E.D.**

**Proof of Theorem 3.** (1) When $t_0^+ \geq t_{right} > 0$ and $\left(\Omega_t^+ X_0, \Phi_{\beta^*}^+ Y_0\right) \in P_{convex}(X,Y)$, $\left(\Omega_t^+ X_0, \Phi_{\beta^*}^+ Y_0\right)$ is located on the weakly efficient frontier $EF_{convex\_weak}$. Assuming $\left(\Omega_t^+ X_0, \Phi_{\beta^*}^+ Y_0\right)$ is not located on the weakly efficient frontier $EF_{convex\_weak}$, we explore the following Model (A.1):

$$\max_{\lambda_j, \theta} \theta$$

$$s.t. \begin{cases} \sum_{j=1}^n \lambda_j x_{ij} = \left(1 + \omega_i t_{right}\right) x_{i0}, i = 1, \dots, m \\ \sum_{j=1}^n \lambda_j y_{rj} \geq \theta \left(1 + \delta_r \beta^*\right) y_{r0}, r = 1, \dots, s \\ \sum_{j=1}^n \lambda_j = 1 \\ \theta \geq 0; \lambda_j \geq 0, j = 1, \dots, n \end{cases} \quad \text{(A.1)}$$

Let $\theta^*$ be the optimal objective of the above model. Because $\left(\Omega_t^+ X_0, \Phi_{\beta^*}^+ Y_0\right)$ is not located on the weakly efficient frontier $EF_{convex\_weak}$, we have $\theta^* > 1$. Therefore, we have $\left(\Omega_t^+ X_0, \Theta^* \Phi_{\beta^*}^+ Y_0\right) \in P_{convex}(X,Y)$, where

$$\Theta^* = diag \underbrace{\left\{\theta^*, \dots, \theta^*\right\}}_{s}.$$

In addition, because $\theta^* > 1$, we have

$$\theta^* \left(1 + \delta_r \beta^*\right) y_{r0} \geq \left(1 + \delta_r \left(\theta^* \beta^*\right)\right) y_{r0}, r = 1, \dots, s$$

Hence, $\theta^* \beta^*$ is a feasible solution for Model (11) and $\theta^* \beta^* > \beta^*$, which contradicts the fact that $\beta^*$ is the optimal solution of Model (11).

(2) When $t_{right} \to 0$, $\left(\Omega_t^+ X_0, \Phi_{\beta^*}^+ Y_0\right)$ converges to $\left(X_0, Y_0\right)$. Thus, there exists a small quantity $t_{right}$ satisfying that both $\left(X_0, Y_0\right)$ and $\left(\Omega_t^+ X_0, \Phi_{\beta^*}^+ Y_0\right)$ obviously have the same supporting hyperplane. This supporting hyperplane may be different for different $t_{right}$. **Q.E.D.**

**Proof of Theorem 4.** We let $\left(\lambda_j^*, \beta^*\right)$ be the optimal solutions of Model (11). Because $t_{right} > 0$ is a sufficiently small quantity, we know that both $\left(X_0, Y_0\right)$ and $\left(\Omega_t^+ X_0, \Phi_{\beta^*}^+ Y_0\right)$ are located on the weakly efficient frontier $EF_{convex\_weak}$, and they have the same supporting hyperplane, where $\Omega_t^+ = diag\left\{1 + \omega_1 t_{right}, ..., 1 + \omega_m t_{right}\right\}$ and $\Phi_{\beta^*}^+ = diag\left\{1 + \delta_1 \beta^*, ..., 1 + \delta_s \beta^*\right\}$. In this case, we know

that the optimal objective value of Model (13) satisfies $\varphi^* = 1$. Thus, we have

$$\begin{cases} U^{T*}Y_0 - V^{T*}X_0 + \mu_0^* = 0 \\ U^{T*}\Phi_{\beta^*}^+ Y_0 - V^{T*}\Omega_t^+ X_0 + \mu_0^* = 0 \end{cases} \qquad (A.2)$$

where $U^* = \left(u_1^*, u_2^*, ..., u_s^*\right)^T$ and $V^* = \left(v_1^*, v_2^*, ..., v_m^*\right)^T$ are the optimal solutions of Model (13).

According to Equation (A.2), we know

$$\sum_{r=1}^{s} u_r^* \delta_r \beta^* y_{r0} - \sum_{i=1}^{m} v_i^* \omega_i t_{right} x_{i0} = 0 \qquad (A.3)$$

From Equation (A.3), we have

$$\frac{\beta^*}{t_{right}} = \frac{\sum_{i=1}^{m} v_i^* \omega_i x_{i0}}{\sum_{r=1}^{s} u_r^* \delta_r y_{r0}} \qquad (A.4)$$

Therefore, when $t_{right} > 0$ is a sufficiently small quantity, we can obtain the optimal objective value of

Model (11) as follows:

$$\xi^* = \frac{\beta^*}{t_{right}} = \frac{\sum_{i=1}^{m} v_i^* \omega_i x_{i0}}{\sum_{r=1}^{s} u_r^* \delta_r y_{r0}} \qquad (A.5)$$

where $U^* = \left(u_1^*, u_2^*, ..., u_s^*\right)^T$ and $V^* = \left(v_1^*, v_2^*, ..., v_m^*\right)^T$ are the optimal solutions of Model (13), and

$\left(U^*, -V^*\right)$ is the normal vector of the supporting hyperplane on the DMU $\left(X_0, Y_0\right)$ and DMU

$\left(\Omega_t^+ X_0, \Phi_{\beta^*}^+ Y_0\right)$.

We know when $0 < t_{right} \leq t_0^+$ and $t_{right} \to 0^+$, both $\left(X_0, Y_0\right)$ and $\left(\Omega_t^+ X_0, \Phi_{\beta^*}^+ Y_0\right)$ are located on the same

"Face" of the weakly efficient frontier, or the value of the Equation (A.5) remains unchanged. Thus, from

Equation (A.5), we know that the optimal objective value of Model (11) is constant with respect to $t_{right}$

when $t_{right} > 0$ is a sufficiently small quantity. **Q.E.D.**

**Proof of Theorem 5.** According to **Definition 9**, we can easily observe that **Theorem 5** holds. **Q.E.D.**

**Proof of Theorem 6.** The proof is similar to that of **Theorem 3** and omitted here.

**Proof of Theorem 7.** The proof is similar to that of **Theorem 4** and omitted here.

**Proof of Theorem 8.** We first discuss the equality between the optimal objective value $\xi^*\left(X_0, Y_0\right)$ and the lower bound $\underline{\rho}\left(X_0, Y_0\right)$. From Equation (A.5), we know that the directional congestion to the "right" of DMU $\left(X_0, Y_0\right)$ reads

$$\xi^*\left(X_0, Y_0\right) = \frac{\beta^*}{t_{right}} = \frac{\sum_{i=1}^{m} v_i^* \omega_i x_{i0}}{\sum_{r=1}^{s} u_r^* \delta_r y_{r0}} \tag{A.6}$$

where $t_{right} > 0$ is a sufficiently small quantity and $\left\{U^* = \left(u_1^*, u_2^*, ..., u_s^*\right)^T, V^* = \left(v_1^*, v_2^*, ..., v_m^*\right)^T, \mu_0^*\right\}$ is the optimal solution of Model (13).

(1)    If    $\xi^*\left(X_0, Y_0\right) < \underline{\rho}\left(X_0, Y_0\right)$    ,    the    optimal    solution    of    Model    (13) $\left\{U^* = \left(u_1^*, u_2^*, ..., u_s^*\right)^T, V^* = \left(v_1^*, v_2^*, ..., v_m^*\right)^T, \mu_0^*\right\}$ satisfies

$$\begin{cases} \sum_{r=1}^{s} u_r^* y_{rj} - \sum_{i=1}^{m} v_i^* x_{ij} + \mu_0^* \leq 0, j = 1, ..., n \\ \sum_{r=1}^{s} u_r^* y_{r0} - \sum_{i=1}^{m} v_i^* x_{i0} + \mu_0^* = 0 \\ \sum_{i=1}^{m} v_i^* x_{i0} - \mu_0^* = 1 \end{cases} \tag{A.7}$$

which contradicts the fact that $\underline{\rho}\left(X_0, Y_0\right)$ is the lower bound in Model (21).

(2) If $\xi^*\left(X_0, Y_0\right) > \underline{\rho}\left(X_0, Y_0\right)$, we can deduce the following formula (A.8) from Equation (A.5):

$$\xi^*\left(X_0, Y_0\right) = \frac{\beta^*}{t_{right}} > \underline{\rho}\left(X_0, Y_0\right) = \frac{\beta_\rho^*}{t_{right}} \tag{A.8}$$

As $t_{right} > 0$, we have $\beta^* > \beta_\rho^*$.

From Model (21), we know that DMU $\left(X_0, Y_0\right)$ satisfies

$$U^{*T} Y_0 - V^{*T} X_0 + \mu_0^* = 0 \tag{A.9}$$

$$U_\rho^{*T} Y_0 - V_\rho^{*T} X_0 + \mu_{\rho 0}^* = 0 \tag{A.10}$$

where $\left(U^*, -V^*\right)$ is the normal vector of a certain "Face" of the weakly efficient frontier on the DMU $\left(X_0, Y_0\right)$, and $\left(U_\rho^*, -V_\rho^*\right)$ is the normal vector of a supporting hyperplane on DMU $\left(X_0, Y_0\right)$.

From (A.6) and (A.9), we have

$$U^{T*} \Phi_{\beta^*}^+ Y_0 - V^{T*} \Omega_t^+ X_0 + \mu_0 = 0 \tag{A.11}$$

where $\Omega_t^+ = diag\{1 + \omega_1 t_{right}, ..., 1 + \omega_m t_{right}\}$ and $\Phi_{\beta^*}^+ = diag\{1 + \delta_1\beta^*, ..., 1 + \delta_s\beta^*\}$.

We pick up a point $\left(\Omega_t^+ X_0, \Phi_{\beta_\rho^*}^+ Y_0\right)$ on the supporting hyperplane with the normal vector $\left(U_\rho^*, -V_\rho^*\right)$ on the DMU $\left(X_0, Y_0\right)$. Thus, we obtain

$$U_\rho^{T*}\Phi_{\beta_\rho^*}^+ Y_0 - V_\rho^{T*}\Omega_t^+ X_0 + \mu_{\rho 0}^* = 0 \qquad (A.12)$$

where $\Phi_{\beta_\rho^*}^+ = diag\{1 + \delta_1\beta_\rho^*, ..., 1 + \delta_s\beta_\rho^*\}$. As DMU $\left(X_0, Y_0\right)$ is strongly efficient, there exists at least one set $\left(U_\rho^*, V_\rho^*\right) > \mathbf{0}$ that satisfies $U_\rho^{T*}Y_0 - V_\rho^{T*}X_0 + \mu_{\rho 0}^* = 0$. In this context, we can obtain the point $\left(\Omega_t^+ X_0, \Phi_{\beta_\rho^*}^+ Y_0\right)$.

As $\beta^* > \beta_\rho^*$, we can obtain the following formula from Equation (A.12):

$$U_\rho^{T*}\Phi_{\beta^*}^+ Y_0 - V_\rho^{T*}\Omega_t^+ X_0 + \mu_{\rho 0}^* > 0 \qquad (A.13)$$

We know that $\left(\Omega_t^+ X_0, \Phi_{\beta^*}^+ Y_0\right)$ is on the weakly efficient frontier $EF_{weak}$, and this fact contradicts the supporting hyperplane $U_\rho^{T*}Y - V_\rho^{T*}X + \mu_{\rho 0}^* = 0$ in Model (21).

Similarly, we can prove that the optimal objective value of Model (15) is equal to the upper bound $\bar{\rho}\left(X_0, Y_0\right)$ in Model (21). **Q.E.D.**

**Proof of Theorem 9.** (1) We suppose that the maximal optimal objective function $\bar{\rho}\left(X_0, Y_0\right)$ of Model (21) is unbounded, but strongly efficient $\left(X_0, Y_0\right)$ is not of the directional smallest scale size. According to Model (15), we know that we can find the optimal solutions of Model (15), denoted by $\left(\beta^*, \eta^*, \lambda_j^*\right)$, in which $\eta^* > 0$. Because $\left(X_0, Y_0\right)$ is strongly efficient, we know that $\beta^* > 0$. According to **Theorems 6 and 7**, we know that we can find a sufficiently small positive constant $t_{left}$ that can ensure that the optimal value of Model (15) is constant. According to **Theorem 8**, we know that the upper bound $\bar{\rho}\left(X_0, Y_0\right)$ of Model (21) is equal to the optimal objective value $\psi^*\left(X_0, Y_0\right)$ of Model (15). This fact contradicts the supposition that the maximal optimal objective function $\bar{\rho}\left(X_0, Y_0\right)$ of Model (21) is unbounded. (2) We can use the similar method to prove this part. **Q.E.D.**

**Proof of Theorem 10. (1) Sufficient condition.** If a DMU $\left(X_0, Y_0\right)$ exhibits directional congestion at the left-hand in the diagonal direction ( $\omega_i = 1, i = 1,...,m$ ; $\delta_r = 1, r = 1,...,s$ ), then there exists $t_{left} > 0$ , which satisfies $\psi^* < 0$ and $\theta^* = 1$ , where $\left( \lambda_j^*, \theta^*, \beta^*, \psi^* \right)$ is the optimal solution of the following Model (A.14):

$$\min_{\lambda_j, \beta} \ \psi = \beta / t_{left}$$

$$s.t. \begin{cases} \sum_{j=1}^{n} \lambda_j x_{ij} = \left( 1 - t_{left} \right) x_{i0}, i = 1,...,m \\ \sum_{j=1}^{n} \lambda_j y_{rj} \geq \left( 1 - \beta \right) y_{r0}, r = 1,...,s \\ \sum_{j=1}^{n} \lambda_j = 1, \lambda_j \geq 0, j = 1,...,n \end{cases} \quad \text{(A.14)}$$

We let $\alpha = 1 - t_{left}, \gamma = 1 - \beta^*$ . Because $\psi^* < 0$ , we have $0 < \alpha < 1, \gamma > 1$ . Thus, there exists $\left( \alpha X_0, \gamma Y_0 \right) \in P_{convex} \left( X, Y \right)$ , which satisfies $\tilde{X}_0 = \alpha X_0 \, (0 < \alpha < 1)$ and $\tilde{Y}_0 \geq \gamma Y_0 \ (\gamma > 1)$ .

**(2) Necessary condition.** If there exists $\left( \tilde{X}_0, \tilde{Y}_0 \right) \in P_{convex} \left( X, Y \right)$ , which satisfies $\tilde{X}_0 = \alpha X_0 \, (0 < \alpha < 1)$ and $\tilde{Y}_0 \geq \gamma Y_0 \ (\gamma > 1)$ , then when we let $\alpha = \left( 1 - t_{left} \right), \gamma = 1 - \beta$ , there exists a negative feasible solution, i.e., the optimal objective value is negative. Therefore, we know that DMU $\left( X_0, Y_0 \right)$ exhibits directional congestion at the left-hand side in the diagonal direction ( $\omega_i = 1, i = 1,...,m$ ; $\delta_r = 1, r = 1,...,s$ ). **Q.E.D.**

**Proof of Theorem 11. (1) Sufficient condition.** If strongly efficient DMU $\left( X_0, Y_0 \right) \in P_{convex} \left( X, Y \right)$ exhibits directional congestion at the left-hand side in a certain direction ( $\left( \omega_1^*,...,\omega_m^* \right)^T$ and $\left( \delta_1^*,...,\delta_s^* \right)^T$ , where $\delta_r^* \geq 0, r = 1,...,s$ and $\omega_i^* \geq 0, i = 1,...,m$ are constants and satisfy $\sum_{r=1}^{s} \delta_r^* = s$ and $\sum_{i=1}^{m} \omega_i^* = m$ ), then we know $\psi^* < 0, \beta^* < 0$ and $\left( \left( 1 - \omega_i^* t_{left} \right) x_{i0}, \left( 1 - \delta_r^* \beta \right) y_{r0} \right) \in P_{convex} \left( X, Y \right)$ where $t_{left} > 0$ and $(\psi^*, \beta^*, \lambda_j^*)$ is the optimal solution of the following Model (A.15):

$$\min_{\lambda_j, \beta} \ \psi = \beta / t_{left}$$

$$s.t. \begin{cases} \sum_{j=1}^{n} \lambda_j x_{ij} = \left( 1 - \omega_i^* t_{left} \right) x_{i0}, i = 1,...,m \\ \sum_{j=1}^{n} \lambda_j y_{rj} \geq \left( 1 - \delta_r^* \beta \right) y_{r0}, r = 1,...,s \\ \sum_{j=1}^{n} \lambda_j = 1, \lambda_j \geq 0, j = 1,...,n \end{cases} \quad \text{(A.15)}$$

Thus, we know that there exists an activity $\left( \left( 1 - \omega_i^* t_{left} \right) x_{i0}, \left( 1 - \delta_r^* \beta \right) y_{r0} \right) \in P_{convex} \left( X, Y \right)$ that uses fewer resources in one or more inputs for making more products in one or more outputs, i.e., weak congestion.

**(2) Necessary condition.** If there exists weak congestion on DMU $\left(X_0, Y_0\right) \in P_{convex}\left(X, Y\right)$, then there exists an activity DMU $\left(\tilde{X}, \tilde{Y}\right) \in P_{convex}\left(X, Y\right)$ that uses less resources in one or more inputs for making more products in one or more outputs, where $\tilde{X} = \left(\tilde{x}_i, i = 1, ..., m\right)$ and $\tilde{Y} = \left(\tilde{y}_r, r = 1, ..., s\right)$. Thus, we know $\tilde{X} \le X_0, \tilde{Y} \ge Y_0$. Let $\theta^*$ be the optimal objective value of the following model:

$$\max_{\lambda_j, \theta} \theta$$
$$s.t. \begin{cases} \sum_{j=1}^{n} \lambda_j X_j = \tilde{X}_0, \sum_{j=1}^{n} \lambda_j Y_j \ge \theta \tilde{Y}_0 \\ \sum_{j=1}^{n} \lambda_j = 1, \lambda_j \ge 0, j = 1, ..., n \end{cases}$$

Thus, we know $\theta^* \ge 1$. In other words, strongly efficient DMU $\left(\tilde{X}, \theta \tilde{Y}\right) \in P_{convex}\left(X, Y\right)$ uses less resources in one or more inputs for making more products in one or more outputs than $\left(X_0, Y_0\right) \in P_{convex}\left(X, Y\right)$. We let $\tilde{x}_i = \left(1 - \omega_i^* t_{left}\right) x_{i0}$ and $\theta \tilde{y}_r = \left(1 - \delta_r^* \beta^*\right) y_{r0}$. Thus, we know $\beta^* < 0, \psi^* = \beta^* / t_{left} < 0$ and

$$t_{left} = \frac{1}{m} \sum_{i=1}^{m} \left(1 - \frac{\tilde{x}_i}{x_{i0}}\right), \quad \beta^* = \frac{1}{r} \sum_{r=1}^{s} \left(1 - \frac{\theta \tilde{y}_r}{y_{r0}}\right)$$

Therefore, there exists a strongly efficient DMU $\left(\tilde{X}, \theta \tilde{Y}\right) \in P_{convex}\left(X, Y\right)$ that uses fewer inputs in the direction $\omega_i^* = \left(1 - \tilde{x}_{i0} / x_{i0}\right) / t_{left}$ for making more outputs in the direction $\delta_i^* = \left(1 - \theta \tilde{y}_{r0} / y_{r0}\right) / \beta^*$. This fact indicates that directional congestion occurs on DMU $\left(X_0, Y_0\right) \in P_{convex}\left(X, Y\right)$ in the above input and output directions. **Q.E.D.**

**Proof of Theorem 12.** If directional congestion exists on the left side in a certain direction on DMU $\left(X_0, Y_0\right)$, we know that there exists at least one strongly efficient DMU $\left(\tilde{X}, \tilde{Y}\right) \in P_{convex}\left(X, Y\right)$, which can produce more outputs using fewer inputs than those of DMU $\left(X_0, Y_0\right)$ in this direction. Because $\sum_{r=1}^{s} \delta_r = s, \delta_r \ge 0; \sum_{i=1}^{m} \omega_i = m, \omega_i \ge 0$, we can observe that there exists an activity in $P_{convex}\left(X, Y\right)$ that uses less resources in one or more inputs for making more products in one or more outputs, i.e., weak or strong congestion occurs. This finding completes the proof. **Q.E.D.**

**Proof of Theorem 13.** According to the definitions of congestion and directional congestion, we can easily see this theorem holds. **Q.E.D.**